\documentclass[11pt]{article}
\usepackage{graphicx}
\usepackage{color}
\usepackage{amssymb,latexsym,amsmath}
\usepackage{epstopdf}
\newtheorem{lemma}{Lemma}[section]
\newtheorem{theo}[lemma]{Theorem}
\newtheorem{prop}[lemma]{Proposition}

\newcommand{\proof}{\noindent{\em Proof: }}

\newcommand{\forme}[1]{}

\def\wbull{\hfill\vrule height .9ex width .8ex depth -.1ex}

\DeclareMathOperator{\PG}{PG}
\def\St{{\mathsf{Star}}}
\def\Li{{\mathsf{Line}}}

\begin{document}

\date{\today}
\title{Derivation of Cameron-Liebler line classes}
\author{{\bf Alexander L. Gavrilyuk}\\
School of Mathematical Sciences,\\ 
University of Science and Technology of China, Hefei 230026, Anhui, PR China\\
and\\
Krasovskii Institute of Mathematics and Mechanics,\\ 
Kovalevskaya str., 16, Ekaterinburg 620990, Russia\\
e-mail: alexander.gavriliouk@gmail.com\\
\\
{\bf Ilia Matkin}\\
Faculty of Mathematics, Chelyabinsk State University,\\
Kashirinykh str., 129, Chelyabinsk 454001, Russia\\
e-mail: ilya.matkin@gmail.com\\
\\
{\bf Tim Penttila}\\
Department of Mathematics, Colorado State University,\\
Fort Collins, CO 80523-1874, USA\\
e-mail: penttila@math.colostate.edu}
\maketitle

\begin{abstract}
We construct a new infinite family of Cameron-Liebler line classes in $\PG(3,q)$ 
with parameter $x=\frac{q^2+1}{2}$ for all odd $q$. 
\end{abstract}

\section{Introduction}\label{intro}
Let $\PG(3,q)$ denote the $3$-dimensional projective space over the finite field $\mathbb{F}_q$.
For a set ${\cal L}$ of lines in $\PG(3,q)$, let ${\overline {\cal L}}$
denote the complementary set of lines. A {\it spread} of $\PG(3,q)$ is a set of $q^2+1$ lines 
that partition the set of points.

We say that ${\cal L}$ is a {\it Cameron-Liebler line class} with parameter $x$ in $\PG(3,q)$, 
if there exists a non-negative integer $x$ such that, for every spread $S$ of $\PG(3,q)$, one has:
\begin{equation*}\label{eq-prop1}
|S\cap \mathcal{L}|=x.
\end{equation*}

It can be seen from the definition that $\overline{\cal L}$ is then a Cameron-Liebler line 
class with parameter $q^2+1-x$, so that we may assume $x\leq \frac{q^2+1}{2}$.
An empty set of lines ($x=0$), the set of all lines in a plane ($x=1$) or, dually, through a point ($x=1$) 
are trivial examples of Cameron-Liebler line classes. If the point is not in the plane, then the union 
of the previous two examples with $x=1$ gives a slightly less trivial Cameron-Liebler line class with parameter $x=2$.

Cameron-Liebler line classes first appeared in the study \cite{CameronLiebler} 
of collineation groups of ${\rm PG}(n,q)$, $n\geq 3$, that have equally many orbits on lines and on points 
(and were given their name in \cite{Penttila}). 
Under the Klein correspondence, Cameron-Liebler line classes are translated to tight sets 
of the Klein quadric being, thus, a special case of a tight set of a polar space 
(see \cite{DrudgeThesis,Metsch-and-Co}).
For more comprehensive background on this topic, we refer to the recent papers \cite{FMX}, 
\cite{Gavrilyuk}, \cite{Metsch2}, \cite{GavrilyukMetsch}, \cite{Metsch-and-Co}.

It was conjectured in \cite{CameronLiebler} that the only Cameron-Liebler line classes are 
the examples mentioned above, i.e., $x\leq 2$. The first counterexample was found by 
Drudge \cite{Drudge} in ${\rm PG}(3,3)$ with $x=5$, which was generalised later 
by Bruen and Drudge \cite{BruenDrudge} 
to an infinite family having parameter $x=\frac{q^2 + 1}{2}$ for all odd $q$. 
The first counterexample in characteristic $2$ was found in \cite{GovaertsPenttila}.
With the aid of computer and using some clever ideas about possible symmetries of Cameron-Liebler 
line classes, Rodgers \cite{Rodgers} constructed many more new examples for certain $x$ and prime powers $q$. 
Very recently, some of them have been shown in \cite{Metsch-and-Co}, \cite{FMX} to be a part of 
a new infinite family of Cameron-Liebler line classes with parameter $x=\frac{q^2+1}{2}$ 
for $q\equiv 5{\rm~or~}9~({\rm mod~}12)$. (In fact, a line class of the family found 
in \cite{Metsch-and-Co}, \cite{FMX} has parameter $\frac{q^2-1}{2}$, however, it is disjoint 
with a plane, which is a Cameron-Liebler line class with parameter $1$, so that the union of their 
lines is a Cameron-Liebler line class with parameter $\frac{q^2-1}{2}+1$.)

In this note, we first describe a switching-like operation in Cameron-Liebler line classes 
that satisfy some necessary conditions (see Lemma \ref{lemma-switching}). We then show 
in Lemma \ref{lemma-x} that these conditions may only hold for line classes with $x=q^2$ 
or $x=\frac{q^2+1}{2}$. Applying this switching operation to the line classes found by Bruen and Drudge, 
we construct another infinite family of Cameron-Liebler line classes 
in ${\rm PG}(3,q)$ with parameter $x=\frac{q^2+1}{2}$ for all odd $q$, and show 
that they are not equivalent to the line classes of Bruen and Drudge, unless $q=3$ 
(see Theorem \ref{theo-main}).


\section{Switching in Cameron-Liebler line classes}
For a point $P$ and a plane $\pi$ of $\PG(3,q)$, let $\St(P)$ and $\Li(\pi)$ 
denote the set of all lines on $P$ or in $\pi$, respectively. 

\begin{lemma}\label{lemma-switching}
Let $\mathcal{L}$ be a Cameron-Liebler line class such that there exists 
an incident point-plane pair $(P,\pi)$ satisfying the following conditions: 
\begin{enumerate}
\item[$(1)$] $(\mathsf{Line}(\pi)\setminus \mathsf{Star}(P))\cap \mathcal{L}=\emptyset$,
\item[$(2)$] $\mathsf{Star}(P)\setminus \mathsf{Line}(\pi)\subseteq \mathcal{L}$.
\end{enumerate}
Then
\[
\mathcal{L}':=\mathcal{L}\cup (\mathsf{Line}(\pi)\setminus \mathsf{Star}(P))\setminus (\mathsf{Star}(P)\setminus \mathsf{Line}(\pi))
\] 
is a Cameron-Liebler line class with the same parameter.
\end{lemma}
\proof For any spread $S$ of $\PG(3,q)$ we have that 
$S$ contains either a line of $\St(P)\cap \Li(\pi)$, or 
a line $\ell\in \Li(\pi)\setminus \St(P)$ and a line $m\in \St(P)\setminus \Li(\pi)$.
In the former case, $S\cap \mathcal{L}=S\cap \mathcal{L}'$, while 
in the latter case $S\cap \mathcal{L}'=(S\cap \mathcal{L})\cup \{m\} \setminus \{\ell\}$.
Thus, $|S\cap \mathcal{L}'|=|S\cap \mathcal{L}|$ holds in both cases, 
and so $\mathcal{L}'$ is a Cameron-Liebler line class. \wbull
\medskip

Let ${\cal L}$ be a Cameron-Liebler line class, and $\ell$ a line of $\PG(3,q)$. 
Then $\ell$ lies in $q+1$ planes $\pi_1$, $\dots$, $\pi_{q+1}$ and contains $q+1$ points $P_1$, $\dots$, $P_{q+1}$. 
Define the square matrix $T(\ell)=(t_{ij})$ of size $q+1$ with integer entries given by
\[
t_{ij}:=|((\Li(\pi_i)\cap\St(P_j))\setminus\{\ell\})\cap \mathcal{L}|,~~1\le i,j\le q+1. 
\]

The set consisting of the matrix $T$, and every matrix obtained from this one 
by a permutation of the rows and a permutation of the columns is called the {\it pattern} of $\ell$ 
with respect to ${\cal L}$. We represent each pattern by one of its matrices. 
This concept was introduced in \cite{Gavrilyuk}, 
where the following result has been proved.

\begin{prop}\label{prop-ti2}
Let $\mathcal{L}$ be a Cameron-Liebler line class with parameter $x$, 
let $\ell$ be a line of $\PG(3,q)$, and $T=(t_{ij})$ the pattern of $\ell$.
\renewcommand{\labelenumi}{\rm (\alph{enumi})}
\begin{enumerate}
\item For any $i\in \{1,\ldots,q+1\}$
\[
\sum_{j=1}^{q+1}t_{ij}=|\Li(\pi_i)\cap \mathcal{L}\setminus \{\ell\}|
\text{~~and~~}
\sum_{j=1}^{q+1}t_{ji}=|\St(P_i)\cap \mathcal{L}\setminus \{\ell\}|.
\]
\item For all $k,l\in\{1,\dots,q+1\}$
$$\sum_{i=1}^{q+1}t_{il}+\sum_{j=1}^{q+1}t_{kj}=
\left\{
\begin{matrix}
x+(q+1)t_{kl}&\mbox{if $\ell\notin \mathcal{L}$} &  & \\
x+(q+1)(t_{kl}+1)-2&\mbox{if $\ell\in \mathcal{L}$}.
\end{matrix}
\right.$$
\item $t_{kl}+t_{rs}=t_{ks}+t_{rl}$ for all $k,l,r,s\in\{1,\dots,q+1\}$.
\item
$$
\sum_{i,j=1}^{q+1} t_{ij}^2=
\left\{
\begin{matrix}
x(q+x)&\mbox{if $\ell\notin \mathcal{L}$} &  & \\
q^3+q^2+(x-1)^2+q(x-1)&\mbox{if $\ell\in \mathcal{L}$}.
\end{matrix}
\right.$$
\end{enumerate}
\end{prop}

\begin{lemma}\label{lemma-x}
Let $\mathcal{L}$ be a Cameron-Liebler line class such that there exists 
an incident point-plane pair $(P,\pi)$ satisfying the conditions 
of Lemma \ref{lemma-switching}. Then the parameter $x$ of $\mathcal{L}$ is equal 
to $q^2$ or $\frac{q^2+1}{2}$.
\end{lemma}
\proof Up to taking the complement to a line set and 
the point-plane duality in $\PG(3,q)$, we may assume 
that there exists a line $\ell$ of $\St(P)\cap \Li(\pi)\setminus \mathcal{L}$.
Let $T$ be the pattern of $\ell$ such that, without loss of generality, 
its first row corresponds to $\pi$, and its first column corresponds to $P$.
Then the conditions of Lemma \ref{lemma-switching} imply that 
$t:=t_{11}=|\St(P)\cap \Li(\pi)\cap \mathcal{L}|$, and $t_{1,j}=q$ and $t_{j,1}=0$ 
for all $j\in \{2,\ldots,q+1\}$. By Proposition \ref{prop-ti2}~(c), we see 
that $t_{ij}=q-t_{11}$ for all $i,j\in\{2,\ldots,q+1\}$.

Further, Proposition \ref{prop-ti2}~(b) applied to the first row and column of $T$, 
and 
Proposition \ref{prop-ti2}~(d) applied to the pattern $T$ give the following equations:
\[
\left\{
\begin{matrix}
t+q^2+t=x+t(q+1),\\
t^2+q^3+q^2(q-t)^2=x(q+x),
\end{matrix}
\right.,
\]
which yield $t=0$ and $x=q^2$ (and thus $\mathcal{L}$ is the complement to 
a Cameron-Liebler line class with parameter $1$), or $t=\frac{q+1}{2}$ and $x=\frac{q^2+1}{2}$.
\wbull

\section{Application of switching}

From Lemma \ref{lemma-x} we see that the only non-trivial case, where the switching 
operation of Lemma \ref{lemma-switching} may be applied, is the case $x=\frac{q^2+1}{2}$.
There exist at least two infinite families of Cameron-Liebler line classes with parameter 
$x=\frac{q^2+1}{2}$, 
namely, the first counterexamples to the Cameron-Liebler conjecture constructed by Bruen and Drudge 
in \cite{BruenDrudge} and the line classes recently found in \cite{Metsch-and-Co} and independently in 
\cite{FMX}. Fortunately, the former satisfy the conditions of Lemma \ref{lemma-switching} 
(while the latter do not), 
and applying the switching operation indeed produces a new Cameron-Liebler line class, 
not equivalent to the original one, if $q>3$. In this section we give the necessary details.

First of all, let us recall the construction by Bruen and Drudge. Let $q$ be an odd prime power, 
and $\mathcal{Q}$ an elliptic quadric of $\PG(3,q)$ with the corresponding quadratic form $\mathsf{Q}$. 
The set of $q+1$ tangents $\mathcal{T}_P$
to a point $P\in \mathcal{Q}$ can be divided into two subsets, say $\mathcal{T}_P^1$, 
$\mathcal{T}_P^2$, of size $(q+1)/2$ each, depending on whether a tangent line 
contains a point $P'\ne P$ such that $\mathsf{Q}(P')$ is a square in $\mathbb{F}_q$.
Note if $\mathsf{Q}(P')$ is a square in $\mathbb{F}_q$, then all the points on the tangent $PP'$
satisfy this property, as $\mathsf{Q}(P+cP')=c^2\mathsf{Q}(P')$.

Denote by $\mathcal{T}^i$ the set $\cup_{P\in \mathcal{Q}}\mathcal{T}_{P}^i$, $i\in \{1,2\}$. 
Let $\mathcal{S}$ and $\mathcal{E}$ be the sets of secant and external lines to $\mathcal{Q}$, 
respectively. Then any of 
\[
\mathcal{S}\cup \mathcal{T}^i,~~
\mathcal{E}\cup \mathcal{T}^j,~~
i,j\in \{1,2\}, 
\]
is a Cameron-Liebler line class of parameter $\frac{q^2+1}{2}$.

Since all these line classes are equivalent under the action of ${\rm P\Gamma L}(4,q)$ 
and the polarity induced by $\mathcal{Q}$ (see \cite{DrudgeThesis}), we may choose, 
without loss of generality, $\mathcal{L}$ to be $\mathcal{S}\cup \mathcal{T}^1$. 
For a point $P_1$ of $\mathcal{Q}$ and its tangent plane $\tau_{P_1}$, 
one can see that  
\[(\mathsf{Line}(\tau_{P_1})\setminus \mathsf{Star}(P_1))\subset \mathcal{E}\subset \overline{\mathcal{L}},~~
\mathsf{Star}(P_1)\setminus \mathsf{Line}(\tau_{P_1})\subset \mathcal{S}\subset \mathcal{L},\]
so that $(P_1,\tau_{P_1})$ satisfies the condition of Lemma \ref{lemma-switching}, 
and the line class $\mathcal{L}'$ defined by
\[
\mathcal{L}':=\mathcal{L}\cup (\mathsf{Line}(\tau_{P_1})\setminus \mathsf{Star}(P_1))\setminus 
(\mathsf{Star}(P_1)\setminus \mathsf{Line}(\tau_{P_1}))
\]
is a Cameron-Liebler line class with parameter $\frac{q^2+1}{2}$.

Our aim now is to show that $\mathcal{L}'$ is not equivalent to $\mathcal{L}$ unless $q=3$.
For $q=3$, we can either apply Drudge's classification of Cameron-Liebler line classes in $\PG(3,3)$ 
\cite{Drudge}
that determined that, up to equivalence, there is a unique Cameron-Liebler line class with parameter $5$, 
or it can be checked with the aid of computer that $\mathcal{L}'$ is projectively equivalent to 
$\overline{\mathcal{L}}$ for this value of $q$.
From now on, we assume that $q>3$.

\begin{lemma}\label{lemma-pl}
A plane $\pi$ of $\PG(3,q)$ contains  
$\frac{q+1}{2}$, or $\frac{q(q+1)}{2}$, or $\frac{(q+1)(q+2)}{2}$ lines of $\mathcal{L}$.
\end{lemma}
\proof If $\pi$ is a tangent plane to $\mathcal{Q}$, then $|\Li(\pi)\cap \mathcal{L}|=\frac{q+1}{2}$ by 
the construction of $\mathcal{L}$. Suppose that $\pi$ is a secant plane so that 
$\pi\cap \mathcal{Q}$ is a conic. Under the polarity, say $\rho$, induced by $\mathcal{Q}$, 
every tangent line to the conic in $\pi$ is mapped to a tangent line to $\mathcal{Q}$ on $\rho(\pi)$.
Therefore, all tangent lines to the conic in $\pi$ are either in $\mathcal{T}^1$ or in $\mathcal{T}^2$. 
In the former case, $\pi$ contains ${q+1 \choose 2}+q+1$ lines from $\mathcal{L}$, 
in the latter case $|\Li(\pi)\cap \mathcal{L}|={q+1 \choose 2}$.
\wbull

\begin{lemma}\label{lemma-pt}
A point $P$ of $\PG(3,q)$ is on 
$q^2+\frac{q+1}{2}$, or $\frac{q(q-1)}{2}$, or $\frac{q(q+1)}{2}+1$ lines of $\mathcal{L}$.
\end{lemma}
\proof If $P\in \mathcal{Q}$, then $|\St(P)\cap \mathcal{L}|=\frac{q+1}{2}+q^2$ by 
the construction of $\mathcal{L}$. Suppose that $P\notin \mathcal{Q}$. 
If $P$ is on a tangent line from $\mathcal{T}^i$ for $i\in \{1,2\}$, 
then all tangent lines to $\mathcal{Q}$ through $P$ are in $\mathcal{T}^i$.
Let $P'$ be a point of $\mathcal{Q}$ such that $PP'$ is a tangent line to $\mathcal{Q}$, 
and consider all secant planes $\pi_1,\ldots,\pi_q$ containing the line $PP'$.
Recall that every point not on a conic in a projective plane of odd order lies on $0$ or $2$ tangents, 
see \cite{Qvist,Segre}. 
Since $\pi_i\cap \mathcal{Q}$ is a conic, and $PP'$ is a tangent line to the conic, 
we conclude that $P$ lies on $2$ tangents and $\frac{q-1}{2}$ secants 
to $\pi_i\cap \mathcal{Q}$. Thus, $|\St(P)\cap \mathcal{L}|=\frac{q(q-1)}{2}$, 
if $PP'\in \mathcal{T}^2$, or $|\St(P)\cap \mathcal{L}|=\frac{q(q-1)}{2}+q+1$, 
if $PP'\in \mathcal{T}^1$.
\wbull

\begin{theo}\label{theo-main}
The line classes $\mathcal{L}$ and $\mathcal{L}'$ are not equivalent under the action 
of ${\rm P\Gamma L}(4,q)$ or a duality.
\end{theo}
\proof Following the notation from the above, one can see that the plane $\tau_{P_1}$ 
contains $\frac{q+1}{2}+q^2$ lines of $\mathcal{L}'$. Since, for a point $P_2\in \mathcal{Q}$, $P_2\ne P_1$, 
one has $\tau_{P_1}\cap \tau_{P_2}\in \mathcal{E}$, the plane $\tau_{P_2}$ contains 
$\frac{q+1}{2}+1$ lines of $\mathcal{L}'$. 
It now follows from Lemmas \ref{lemma-pl}, \ref{lemma-pt} that 
the intersection numbers of $\mathcal{L}'$ with respect to planes and points of $\PG(3,q)$ 
are different from those of $\mathcal{L}$ or $\overline{\mathcal{L}}$.\wbull
\medskip

We also note that $\mathcal{L}'$ is not equivalent to a line class of the family found in 
\cite{Metsch-and-Co}, \cite{FMX}, since there is no plane (or, dually, a point with all lines on it) contained in or disjoint from $\mathcal{L}'$. 
In particular, in $\PG(3,5)$, there exist at least three pairwise non-equivalent Cameron-Liebler 
line classes with $x=\frac{q^2+1}{2}=13$ (namely, the example by Bruen and Drudge, its switched mate 
by Theorem \ref{theo-main}, and the example found in \cite{FMX} and \cite{Metsch-and-Co}). 
In fact, up to equivalence, these are the only Cameron-Liebler line classes with given $x$ in $\PG(3,5)$ 
(the details will be given elsewhere). 

The line class $\mathcal{L}'$ contains only the one incident point-plane pair, 
namely, $(P_1,\tau_{P_1})$, satisfying the conditions of Lemma \ref{lemma-switching}, and, clearly, 
switching of $\mathcal{L}'$ with respect to it gives the line class $\mathcal{L}$.
Since, for $q>3$, there is a unique switched mate for $\mathcal{L}'$ (namely, $\mathcal{L}$), 
it follows that its stabiliser $G_{\mathcal{L}'}$ is a subgroup of the stabiliser $G_{\mathcal{L}}$. 
The stabiliser $G_{\mathcal{L}}$ of a Bruen-Drudge line class is a subgroup of index two of 
${\rm P\Gamma O}^-(4,q)$, i.e., the subgroup that fixes $\mathcal{T}^1$ and $\mathcal{T}^2$. 
Thus, $G_{\mathcal{L}'}$ is the stabiliser of the point $P_1$ in $G_{\mathcal{L}}$. 
Then, for $q=p^h$, where $p$ is a prime, $G_{\mathcal{L}'}$ has order $q^2(q^2-1)h$, and 
is isomorphic to ${\rm AGL}(1,q^2)\rtimes C_h$.

We expect that the only non-trivial Cameron-Liebler line classes satisfying the conditions 
of Lemma \ref{lemma-switching} are the examples of Bruen and Drudge and their switched mates.

\noindent{\bf Acknowledgments.}
The research of A.L.G. was funded by Chinese Academy of Sciences President's International 
Fellowship Initiative (Grant No. 2016PE040). I.M. was partially supported by the Russian 
Foundation for Basic Research (Grant No. 16-31-50070).

The Cameron-Liebler line classes in this paper were constructed independently at least three times: 
first by Penttila who didn't publish the result, but did tell others of the construction, then by 
Cossidente and Pavese in \cite{Cossidente}, and finally by Gavrilyuk and Matkin. Gavrilyuk recalled 
Penttila mentioning something to him about the construction and contacted him with the outcome being 
the decision to write up the construction jointly. Since the approach in this paper shows a different 
point of view to that in the paper by Cossidente and Pavese, we felt that the construction 
in this paper deserves to be published.
\smallskip

We would like to thank Anton Betten for organising the Combinatorics and Computer Algebra 2015 
conference, whose problem sessions brought us together to begin this work. 
Part of the work was done while I.M. was visiting USTC, he thanks Jack Koolen for his hospitality. 


\begin{thebibliography}{}
\bibitem{Metsch-and-Co}
J.~De~Beule,~J.~Demeyer,~K.~Metsch,~M.~Rodgers.
\newblock A new family of tight sets in $Q^+(5,q)$. 
\newblock {\em Des. Codes Cryptogr.} \textbf{78} (2016) 655--678.

\bibitem{BruenDrudge}
A.~A. Bruen,~Keldon~Drudge.
\newblock The construction of {C}ameron-{L}iebler line classes in {${\rm PG}(3,q)$}.
\newblock {\em Finite Fields Appl.} \textbf{5(1)} (1999) 35--45.

\bibitem{CameronLiebler}
P.~J. Cameron,~R.~A. Liebler.
\newblock Tactical decompositions and orbits of projective groups.
\newblock {\em Linear Algebra Appl.} \textbf{46} (1982) 91--102.

\bibitem{Cossidente}
A.~Cossidente,~F.~Pavese.
\newblock Intriguing sets of quadrics in $\PG(5,q)$.
\newblock {\em Adv. Geom.}, in press.


\bibitem{DrudgeThesis} 
Keldon Drudge.
\newblock Extremal sets in projective and polar spaces.
\newblock {\em Ph.D. Thesis}, University of Western Ontario, 1998.

\bibitem{Drudge}
Keldon Drudge.
\newblock On a conjecture of {C}ameron and {L}iebler.
\newblock {\em European J. Combin.} \textbf{20(4)} (1999) 263--269.

\bibitem{FMX}
T.~Feng,~K.~Momihara,~Q.~Xiang.
\newblock {C}ameron-{L}iebler line classes with parameter $x=\frac{q^2-1}{2}$. 
\newblock {\em J. Combin. Theory Ser. A} \textbf{133} (2015) 307--338.

\bibitem{GavrilyukMetsch}
A.~L.~Gavrilyuk,~K.~Metsch.
\newblock A modular equality for {C}ameron-{L}ieber line classes. 
\emph{J. Combin. Theory Ser. A} \textbf{127} (2014) 224--242.

\bibitem{Gavrilyuk}
Alexander~L. Gavrilyuk, ~Ivan~Yu.~Mogilnykh.
\newblock Cameron-{L}iebler line classes in {${\rm PG}(n,4)$}.
\newblock {\em Des. Codes Cryptogr.} \textbf{73(3)} (2014) 969--982.

\bibitem{GovaertsPenttila}
Patrick~Govaerts,~Tim~Penttila.
\newblock Cameron-{L}iebler line classes in {${\rm PG}(3,4)$}.
\newblock {\em Bull. Belg. Math. Soc. Simon Stevin} \textbf{12(5)} (2005) 793--804.





\bibitem{Metsch2}
Klaus Metsch.
\newblock An improved bound on the existence of {C}ameron--{L}iebler line classes.
\newblock {\em J. Combin. Theory Ser. A} \textbf{121} (2014) 89--93.

\bibitem{Penttila}
Tim Penttila.
\newblock Cameron-{L}iebler line classes in {${\rm PG}(3,q)$}.
\newblock {\em Geom. Dedicata} \textbf{37(3)} (1991) 245--252.

\bibitem{Qvist}
Bertil Qvist.
\newblock Some remarks concerning curves of the second degree in a finite plane.
\newblock {\em Suomalainen tiedeakatemia (Sci. Fennica)} \textbf{134} (1952) 4--27.


\bibitem{Rodgers}
Morgan Rodgers.
\newblock Cameron-{L}iebler line classes.
\newblock {\em Des. Codes Cryptogr.} \textbf{68(1-3)} (2013) 33--37.

\bibitem{Segre}
Beniamino Segre.
\newblock Ovals in a finite projective plane.
\newblock {\em Canad. J. Math.} \textbf{7} (1955) 414--416.
\end{thebibliography}


\end{document}